\documentclass[3p]{elsarticle}

\usepackage{amssymb}
\usepackage{amsthm}

\usepackage{amsmath,amsthm,a4wide,mathrsfs,color,array}
\usepackage[utf8]{inputenc}
\usepackage[T1]{fontenc}
\usepackage{subcaption}
\usepackage{graphicx}
\usepackage{float}
\usepackage{multirow}
\usepackage{hyperref}
\usepackage{cleveref}
\usepackage{tikz-cd}
\usepackage{fancyvrb}
\usepackage{xcolor}
\usetikzlibrary{shapes,arrows,positioning}

\usepackage{tikz}
\usetikzlibrary{shapes,arrows,positioning}
\usetikzlibrary{babel}
\usetikzlibrary{shapes,arrows}
\definecolor{dark_blue_pers}{RGB}{46,87,144}
\definecolor{blue_pers}{RGB}{54,104,171}
\definecolor{grey_pers}{RGB}{245,245,245}
\definecolor{red_pers}{RGB}{213,78,33}
\usetikzlibrary{backgrounds, scopes}

\usepackage{color}
\definecolor{dark_blue}{RGB}{46,87,144}
\definecolor{dark_green}{RGB}{0,100,0}

\newcommand{\rev}{}

\newcommand{\sigmoid}{\operatorname{sigmoid}}
\newcommand{\argmin}{\operatorname{argmin}}
\newcommand{\nc}{\newcommand}

\nc{\norm}[2]{\left\|#1\right\|_{#2}}

\usepackage{algorithm}
\usepackage{algorithmic}

\nc{\IC}{\mathbb{C}}
\nc{\IE}{\mathbb{E}}
\nc{\IN}{\mathbb{N}}
\nc{\IR}{\mathbb{R}}

\nc{\be}{\begin{equation}}
\nc{\ee}{\end{equation}}

\nc{\mN}{\mathcal{N}}
\nc{\mB}{\mathcal{B}}
\nc{\mT}{\mathcal{T}}
\nc{\mL}{\mathcal{L}}
\nc{\mE}{\mathcal{E}}

\nc{\bx}{{\bf x}}
\nc{\bn}{{\bf n}}
\nc{\bq}{{\bf q}}
\nc{\bW}{{\bf W}}
\nc{\bb}{{\bf b}}
\nc{\bz}{{\bf z}}
\nc{\bu}{{\bf u}}
\nc{\btheta}{{\boldsymbol{\theta}}}

\nc{\loc}{{_\textup{loc}}}

\usepackage[normalem]{ulem}

\journal{}

\begin{document}

\begin{frontmatter}
%Efficient?

\title{Hyper-parameter tuning of physics-informed neural networks:\\Application to Helmholtz problems}

\author[rvt,do]{Paul Escapil-Inchausp\'e\corref{cor1}}
\ead{paul.escapil@edu.uai.cl}
\author[rvt,do,capes]{Gonzalo A. Ruz}
\ead{gonzalo.ruz@uai.cl}

\address[rvt]{Facultad de Ingenier\'ia y Ciencias, Universidad Adolfo Ib\'a\~nez, Av.~Diagonal Las Torres 2640, Pe\~nalol\'en, Santiago, Chile}
\address[do]{Data Observatory Foundation, Santiago, Chile}
\address[capes]{Center of Applied Ecology and Sustainability (CAPES), Santiago, Chile}
\cortext[cor1]{Corresponding author}

\begin{abstract}
We consider physics-informed neural networks \rev{(PINNs)} [Raissi et al., \emph{J.~Comput.~Phys.} 278 (2019) 686-707] for forward physical problems. In order to find optimal PINNs configuration, we introduce a \emph{hyper-parameter \rev{optimization}} \rev{(HPO)} procedure via Gaussian processes-based Bayesian optimization. We apply the \rev{HPO} to Helmholtz \rev{equation} for bounded domains and conduct a thorough study, focusing on: (i) performance, (ii) the collocation points density $r$ and (iii) the frequency $\kappa$, confirming the applicability and necessity of the method. Numerical experiments are performed in two and three dimensions, including comparison to finite element methods.
\end{abstract}


\end{frontmatter}
%\tableofcontents 
\section{Introduction}\label{sec:intro}
Having efficient methods to solve physical problems is key in several fields ranging from electronic design automation to optics and acoustics \cite{Nedelec}. Physics can often be described by partial differential equations (PDEs) with suitable boundary conditions (BCs), i.e.~as boundary value problems (BVPs) \cite{Nedelec,steinbach2007numerical}. Under appropriate conditions on the domain and the source term, BVPs are known to be well-posed on a continuous level \cite{ern2004theory}. Amongst other physical problems, acoustic wave behavior is often described by Helmholtz equations \cite{Nedelec}\rev{, whose underlined operator is coercive \cite[Section 3.6]{steinbach2007numerical}---of the form elliptic+compact operator.} 

Traditional schemes for solving BVPs include finite element methods (FEM)   \cite{steinbach2007numerical,ern2004theory}, spectral methods or boundary element methods (BEM) \cite{chandler2012numerical,sauter}, the latter being commonly used for unbounded domains.

These techniques benefit from an enriched theory, including precise convergence bounds for both the solution error and iterative solvers \cite{steinbach2007numerical}. \rev{They have been the state-of-the art solution in engineering applications over the past decades}. However, their resolution can be numerically expensive (in particular for growing wavenumbers in propagation problems \cite{babuskapollution2000}), and they do not naturally incorporate additional information (e.g., sensors). Furthermore, they are linear by nature (\rev{i.e.~}they involve bounded linear operators), and they do not adapt with high dimension\rev{al} problems.

As opposed to the previous techniques, deep learning (DL) \cite{bengio2017deep} has shown to be a promising research area. Recent advances in hardware capabilities (GPU acceleration, increase in computational power) made it possible to use deep neural networks (DNNs) to represent complex problems. They allow simulating high-dimensional mappings \cite{lu2021deepxde,khoo_lu_ying_2021} with direct application to uncertainty quantification \cite{scarabosio2021deep}. 

One type of DNNs are physics-informed neural networks which were introduced recently in \cite{RAISSI2019686}. They encode the boundary value problem in the loss function and rely on automatic differentiation. Their strength is that they combine the \rev{aforementioned strengths of both DL and classical numerical analysis.}%generality of DL with the specificity of classical numerical analysis. 

The physics (PDE, BCs), and more generally, the additional information the user has about the problem, can be plugged into the DNN. This allows for genericity and simplicity. Among others, PINNs have proved to be useful to solve inverse problems \cite{Chen2020} or stochastic PDEs \cite{chen2021learning,ZHANG2019108850,MENG2020109020}. Also refer to \cite{Wang2020WhenAW} and the references therein.

Several active open-source libraries have been created, demonstrating a growing interest in PINNs. We mention DeepXDE \cite{lu2021deepxde} and \cite{zubov2021neuralpde,mcclenny2021tensordiffeq,haghighat2021sciann}. These are built on top of DL libraries such as TensorFlow, Keras, or PyTorch. 

By nature, PINNs exhibit a number of \emph{hyper-parameters} (HPs) such as: the learning rate, the width and depth for the DNN, the activation function and the weights for the loss function. Notice that one could include more HPs, such as the optimization algorithm, the learning rate function complexity (e.g.~using a learning rate scheduler), or the collocation points distribution. The high dimensionality for the HPs search space makes it difficult to find proper configurations\rev{---this means configurations that lead to satisfactory generalization results.} \rev{Furthermore, slow training has been a frequent concern for the PINNs community \cite{haghighat2021sciann,markidis2021old}.}

\rev{To remedy those concerns, HPO} \cite{NIPS2011_86e8f7ab}---also referred to as (HP) tuning---is an insightful solution. Many approaches exist, such as grid (resp.~random) search optimization \cite{NIPS2011_86e8f7ab}. Still, \rev{grid search suffers} from the curse of dimensionality, and \rev{random search} can be very inefficient. \rev{Furthermore, each assessment of HPs leads to train a DNN, advocating for limiting the total number of configurations. To mitigate these drawbacks, (Gaussian processes-based) Bayesian optimization \cite{snoek2012practical,yu2020hyper} is a method of choice. It consists in applying GP-based regression to the previous assessed HPs in order to predict the next (most likely) best configuration.}

Finally, it is known that ``the direct usage of PINN in scientific applications is still far from meeting computational performance
and accuracy requirements'' \cite{Chen2020}, as compared to traditional solvers. It could be studied more accurately how far GPU-PINN is from meeting traditional FEM (e.g., in terms of orders of magnitude, for both the memory requirements and computational times).

In this work, we apply HPO via Gaussian processes-based Bayesian optimization to mitigate the poor training of PINNs. We focus on forward problems for Helmholtz operator. \rev{It is chosen as it is: a linear (yet coercive) operator, which encompasses the issues related to the oscillatory modes \cite{babuskapollution2000}; and it paves the way towards more complex cases including fluid mechanics and electromagnetic waves simulation \cite{Nedelec}.}

Also, to our knowledge, application of PINNs to this operator was not studied thoroughly before (except from \cite{Chen2020,luluhard2021}). We conduct an exhaustive overview of tuning concerning: (i) performance, (ii) the collocation point density ---or precision---$r$ and (iii) the wavenumber $\kappa$. We consider two- (resp.~three-) dimensional cases with Dirichlet (resp.~Neumann) BCs, and compare the PINNs to FEM. Numerical results are performed with the state-of-the-art DeepXDE\footnote{\url{https://github.com/lululxvi/deepxde/}} Python library \cite{lu2021deepxde}. \rev{Remark that a tutorial was added by the authors to DeepXDE documentation.}

\rev{Our main practical findings are as follows:} HPO leads to a significant reduction for the loss (e.g.~$7$ orders of magnitude in two dimensions for $\kappa =4\pi$), hence the necessity of the tuning procedure. According to HPO for Helmholtz problems, our numerical results issue the following practical guidance for the PINNs configuration:
\begin{itemize}
    \item Quite shallow DNNs, with depth $L-1 \in \{2,3\}$, and constant width $N \in [250,500]$;
    \item Learning rates $\alpha \in [10^{-4}, 10^{-3}]$, and $\sigma= \sin$ activation function;
    \item Small BCs loss term $w_\Gamma \in [1, 10]$.
\end{itemize}
Moreover, we observe that the training deteriorates with increasing $\kappa$---due to the frequency-principle \cite{markidis2021old,xu2019frequency}---and precision $r$. 

%With respect to the performance of the tuned PINNs, we highlight %the following:
%\begin{itemize}
%    \item They can lead to small residuals with few iterations, %e.g.~around $5,000$;
%    \item They can require few collocation points with respect %to the wavelength, e.g.~can lead to $20\%$ relative norm error %with only $5$ collocation points per wavelength;
%    \item They prove to scale well with the dimension $d$, and %can compete with classic FEM for $d=2,3$.
%\end{itemize}

This work is structured as follows: we introduce PINNs in \Cref{sec:PINNs} and the Bayesian HPO in \Cref{sec:HPO}. Consequently, we conduct exhaustive numerical experiments in \Cref{sec:NumExp} and discuss further research avenues in \Cref{sec:Conclusion}.

\section{Physics-informed neural networks (PINNs)}\label{sec:PINNs}
Set $d\geq 1$ and consider a bounded computational domain $D\subset \IR^d$ with boundary $\Gamma := \partial D$ and exterior unit normal field $\bn$. 

For any scalar field $u$, consider the following problem:
\begin{align*}
\mN[u] = f \quad \text{in}\quad D,\\
\mB[u] = g \quad \text{on} \quad \Gamma,
\end{align*}
with $\mN$ a potentially nonlinear differential operator, and $\mB$ the boundary conditions (BCs) operator. For example, the Helmholtz operator reads:
$$
\mN[u] := - \Delta u - \kappa^2 u,
$$
with $\kappa$ the wavenumber. As BCs, we introduce:
\be
\begin{array}{lll}
\text{(Dirichlet BC): }& \mB[u] := u|_\Gamma = g \quad \text{on}\quad \Gamma,\\
\text{(Neumann BC): } & \mB[u]:=\nabla u|_\Gamma \cdot \bn = g \quad \text{on}\quad \Gamma.
\end{array}
\ee
\rev{Throughout this manuscript, we assume that that the problems under consideration are well-posed (refer to \cite[Section 4.4]{steinbach2007numerical} for more details)}. Let $\sigma $ be a smooth activation function. Given an input $\bx \in D \subset \IR^d$ and following \cite[Section 2.1]{lu2021deepxde}, we define an $L$-layer neural network with $N_l$ neurons in the $l$-th layer for $1 \leq l \leq L-1$ ($N_0= d$ and $N_L = 1$). For $1 \leq l \leq L$, let us denote the weight matrix and bias vector in the $l$-th layer by $\bW^l \in \IR^{N_l \times N_{l-1}}$ and $\bb^l \in \IR^{N_l}$, respectively. The solution $u$ can be approximated by a deep (feedforward) neural network defined as follows:
\be 
\begin{array}{rll}
\text{input layer:} \quad & \bx \in \IR^d,\\
\text{hidden layers:} \quad & \bz^l (\bx) = \sigma( \bW^l \bz^{l-1} (\bx) + \bb^l) \in \IR^{N_l},\\
 & \quad \text{for} \quad 1 \leq l \leq L-1,\\
\text{output layer:} \quad & \bz^L(\bx) = \bW^{\rev{L}} \bz^{L-1} (\bx) + \bb^{\rev{L}}  \in \IR. 
\end{array}
\ee
This results in the representation $u_\theta(\bx) : = \bz^L(\bx)$, with 
$$\theta:= \{ (\bW^1, \bb^1), \cdots, (\bW^L, \bb^L)\}
$$
the (trainable) parameters---or weights---in the network. For the sake of simplicity, we set:
\be
\Theta := \IR^{|\theta|}.
\ee
 Acknowledge that \cite[Eq. (2.11)]{Mishra2020EstimatesOT}
\be\label{eq:trainable}
|\Theta| =: \dim(\Theta)= \sum_{l = 1}^{L} N_{l} (N_{l-1} + 1).
\ee
Next, we introduce collocation points (sensors) for both the domain $\mT_D:= \{\bx_i^D\}_{i=1}^{|\mT_D|}$ and the boundary $\mT_\Gamma :=\{\bx_i^\Gamma\}_{i=1}^{|\mT_\Gamma|}$, and assume that we have observations $\{u_i(\bx^u_i)\}_{i=1}^{|\mT_u|}$ on $\mT_u := \{\bx^u_i\}_{i=1}^{|\mT_u|}$. Finally, we set 
$$\mT: = \mT_D \cup \mT_\Gamma \cup \mT_u.$$
Notice that the formulation remains valid for unlabelled data, i.e.~$\mT_u = \{\}$.

\begin{figure}[!htb]
\centering
\resizebox{8.5cm}{!} {
 \begin{tikzpicture}[node distance = 2cm, auto]
   % Draw the NN
   \node[rectangle, draw, minimum size = 5cm, rounded corners = .8ex, dashed] (boxNN)   at (-8,0){};
   \node[rectangle, draw, minimum size = 5cm, rounded corners = .8ex, dashed] (boxPDE)   at (-2,0){};
    \node[circle, draw, minimum size = .5cm, fill = black!15] (x)   at (-10,0)  {$x$};
    \node[circle, draw, minimum size = .5cm, fill = white] (S1)  at (-8.5,-1.5) {$\sigma$};
    \node[circle, draw, minimum size = .5cm, fill = white] (S2)  at (-8.5,-.75) {$\sigma$};
    \node[circle, draw, minimum size = .5cm, fill = white] (S3)  at (-8.5,0) {$\sigma$};
    \node[circle, draw, minimum size = .5cm, fill = white] (S4)  at (-8.5,.75) {$\sigma$};
    \node[circle, draw, minimum size = .5cm, fill = white] (S5)  at (-8.5,1.5) {$\sigma$};
    \node[circle, draw, minimum size = .5cm, fill = white] (S6)  at (-7.5,-1.5) {$\sigma$};
    \node[circle, draw, minimum size = .5cm, fill = white] (S7)  at (-7.5,-.75) {$\sigma$};
    \node[circle, draw, minimum size = .5cm, fill = white] (S8)  at (-7.5,0) {$\sigma$};
    \node[circle, draw, minimum size = .5cm, fill = white] (S9)  at (-7.5,.75) {$\sigma$};
    \node[circle, draw, minimum size = .5cm, fill = white] (S10) at (-7.5,1.5) {$\sigma$};
    \node[circle, draw, minimum size = .5cm, fill = black!15] (u)   at (-6,0) {$u_\theta$};

   \begin{scope}[on background layer]  
    \draw [-stealth] (x) -- (S1);
    \draw [-stealth] (x) -- (S2);
    \draw [-stealth] (x) -- (S3);
    \draw [-stealth] (x) -- (S4);
    \draw [-stealth] (x) -- (S5);

    \draw [-stealth] (S1) -- (S6);
    \draw [-stealth] (S1) -- (S7);
    \draw [-stealth] (S1) -- (S8);
    \draw [-stealth] (S1) -- (S9);
    \draw [-stealth] (S1) -- (S10);

    \draw [-stealth] (S2) -- (S6);
    \draw [-stealth] (S2) -- (S7);
    \draw [-stealth] (S2) -- (S8);
    \draw [-stealth] (S2) -- (S9);
    \draw [-stealth] (S2) -- (S10);

    \draw [-stealth] (S3) -- (S6);
    \draw [-stealth] (S3) -- (S7);
    \draw [-stealth] (S3) -- (S8);
    \draw [-stealth] (S3) -- (S9);
    \draw [-stealth] (S3) -- (S10);

    \draw [-stealth] (S4) -- (S6);
    \draw [-stealth] (S4) -- (S7);
    \draw [-stealth] (S4) -- (S8);
    \draw [-stealth] (S4) -- (S9);
    \draw [-stealth] (S4) -- (S10);

    \draw [-stealth] (S5) -- (S6);
    \draw [-stealth] (S5) -- (S7);
    \draw [-stealth] (S5) -- (S8);
    \draw [-stealth] (S5) -- (S9);
    \draw [-stealth] (S5) -- (S10);

    \draw [-stealth] (S6) -- (u);
    \draw [-stealth] (S7) -- (u);
    \draw [-stealth] (S8) -- (u);
    \draw [-stealth] (S9) -- (u);
    \draw [-stealth] (S10) -- (u);
    \end{scope}
    \node[text width = 3cm] at (-6, 2.2) {{\bf NN:} $u_\theta(x)$};

    % Draw the PDE
    \node[circle, draw, minimum size = 1.25cm, fill = white] (L)   at (-3.5,1)  {$\mN[\cdot, k]$};
    \node[circle, draw, minimum size = 1.25cm, fill = white] (B)   at (-3.5,-1)  {$\mB[\cdot]$};

    \node[rectangle, draw, minimum size = 1.25cm, fill = blue_pers!20, minimum width = 2.3cm] (resL)   at (-1,1)  {$\mN[u_\theta, k]-f$};
    \node[rectangle, draw, minimum size = 1.25cm, fill = blue_pers!20, minimum width = 2.3cm] (resB)   at (-1,-1)  {$\mB[u_\theta]-g$};    

    \draw [-stealth] (u) -- (L);
    \draw [-stealth] (u) -- (B);

    \draw [-stealth] (L) -- (resL);
    \draw [-stealth] (B) -- (resB);

    \node[text width = 3cm] at (0, 2.2) {{\bf PDE:}};

    % Loss and min
    \node[rectangle, draw, minimum size = 1cm, fill = red_pers!20, minimum width = 1cm, rounded corners = 0.4ex] (res)   at (1.5,0)  {$\mL_\theta$};    
    \draw [-stealth] (resL) -- (res);
    \draw [-stealth] (resB) -- (res);

    \coordinate (f1) at (1.5, 3);
    \coordinate (f2) at (-8, 3);

    \draw [-] (res) -- (f1);
    \draw [-] (f1) -- (f2);
    \draw [-stealth] (f2) -- (boxNN);
    \node[text width = 3cm] at (-2, 3.25) {{\bf Optimization:} $\theta^\star$};

\end{tikzpicture}}
\caption{Schematic representation of a PINN. A NN with $L=3$ (i.e. $L-1 = 2$ hidden layers) and $N=5$ learns the mapping $\bx \mapsto u(\bx)$. The PDE is taken into account throughout the residual $\mL_\theta$, and the NN trainable parameters $\theta$ are optimized via NN training, leading to optimal $\theta^\star$.}
\label{picture:PINNs}
\end{figure}
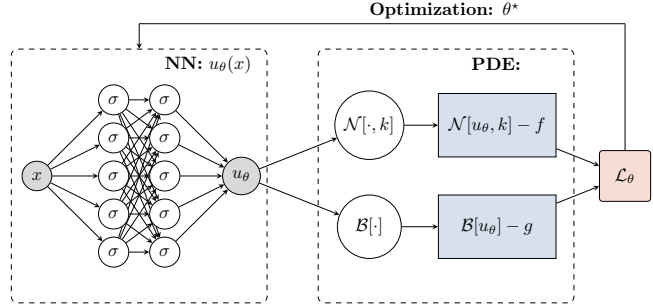

For any $w_D,w_\Gamma,w_u > 0$, the weighted composite loss function is defined as:
\be\label{eq:weightedCompositeLoss}
\mL_\theta = w_D \mL^D_\theta + w_\Gamma\mL^\Gamma_\theta  +w_u \mL^u_\theta,
\ee
where
$$
\mL^D_\theta := \frac{1}{N_D}\sum_{\bx \in \mT_D} \big{|}(\mN[u_\theta,k]- f)(\bx)\big{|}^2,
$$
$$
\mL^\Gamma_\theta := \frac{1}{N_\Gamma}\sum_{\bx \in \mT_\Gamma} \big{|}(\mB[u_\theta] - g) (\bx)\big{|}^2
$$
and
$$
\mL^u_\theta := \frac{1}{N_u}\sum_{\bx \in \mT_u} \big{|}(u_\theta- u_i)(\bx)\big{|}^2.
$$
A schematic representation of a PINN is shown in \Cref{picture:PINNs}.

Remark that the PDE and BCs are incorporated throughout the loss function, and evaluated via automatic differentiation. The optimization (training) procedure allows to define the control weights $\theta\in \Theta$ of the network. For this purpose, the sensors are partitioned into a training and testing data:
$$
\mT_\cdot^\text{train} \quad \text{and} \quad \mT_\cdot^\text{test},
$$
respectively, with $\cdot$ being either $D,\Gamma,u$. Accordingly, we introduce $\mT^\text{train}$ and $\mT^\text{test}$, and the losses $\mL_\theta^\text{train}$ and $\mL_\theta^\text{test}$. We define the target as the training loss function $\mL_\theta^\text{train}$. In parallel, when we are with an exact solution, we use the relative \rev{$l^2$}-norm error taken on the testing sensors $\mT_D^\text{test}$:
$$
\mL_\theta^\text{metric}  = \frac{\|\bu_\theta - \bu\|_2}{\|\bu\|_2},
$$ 
with vectors $\bu,\bu_\theta \in \IR^{|\mT^\text{test}_D|}$ with coefficients $\bu_i:= u(\bx_i^D)$ and $\bu_{\theta,i} = u_\theta(\bx_i^D)$ respectively.

For Dirichlet BCs, one can enforce hard constraint BCs \cite[Section 2.3]{luluhard2021} by applying a transformation to the net:
\be\label{eq:Hard1}
\hat{u}_\theta (\bx):  = g(\bx) + \rev{\ell}  (\bx) u_\theta(\bx)\quad \bx \in \overline{D}
\ee
where 
\be\label{eq:Hard2} 
\begin{array}{rll}
\rev{\ell}(\bx) = 0 , & \bx \in \Gamma, \\
\rev{\ell}(\bx) > 0 , & \bx \in D.
\end{array}
\ee

We are seeking:
\be\label{eq:InnerOptimizationProblem}
\theta^\star = \argmin_{\theta \in \Theta}(\mL_\theta)
\ee
as being the \emph{inner optimization problem}. Next, the HPO will introduce an outer counterpart. An approximation to \eqref{eq:InnerOptimizationProblem} is delivered via an iterative optimizer such as ADAM (see \Cref{alg:ADAM}). Notice that further application of L-BFGS \cite{lbfgs} can improve training \cite{lu2021deepxde}. 
\begin{algorithm}[!ht]
\caption{ADAM \cite[Algorithm 1]{kingma2014adam} applied to loss function $\mL(\theta)$. $g_t^2$ indicates the elementwise square. In this work, we use the default settings for $\beta_1 = 0.9$, $\beta_2 = 0.999$ and $\epsilon = 10^{-7}$ (resp.~$\epsilon = 10^{-8}$) for TensorFlow (resp.~PyTorch).}
\begin{algorithmic}
\label{alg:ADAM}
\REQUIRE $\alpha$ (learning rate)
\REQUIRE $\beta_1,\beta_2 \in [0,1)$ (exponential decay rates for the moment estimates)
\REQUIRE $\theta_0$ (initial trainable parameters) 
\STATE Initialize 1$^{st}$ and 2$^{nd}$ moment vectors: $ m_0 \leftarrow0$, $v_0 \leftarrow 0$
\FOR{$k=0$ \TO $K-1$} 
\STATE $g_{k+1} \leftarrow \nabla l (\theta_k)$
\STATE $m_{k+1} \leftarrow \beta_1 \cdot m_{k} + (1- \beta_1) \cdot g_{k+1}$ 
\STATE $v_{k+1} \leftarrow \beta_2 \cdot v_{k} + (1-\beta_2) \cdot g_{k+1}^2$
\STATE $\hat{m}_{k+1} \leftarrow m_{k+1} / (1 - \beta_1^{k+1})$
\STATE $\hat{v}_{k+1} \leftarrow v_{k+1} / (1 - \beta_2^{k+1})$
\STATE $\theta_{k+1} \leftarrow \theta_{k} - \alpha \cdot \hat{m}_{k+1} / (\sqrt{{\hat{v}_{k+1}}} + \epsilon)$
\ENDFOR
\ENSURE $\theta_K^+ \leftarrow \argmin_{ k \in \{0 ,\cdots, K \}} l(\theta_k)$
\end{algorithmic}
\end{algorithm}

Application of ADAM optimizer to the PINNs with loss $\mL_\theta^\text{train}$ leads to:
\be 
\theta_K^+ := \argmin_{ k \in \{0 ,\cdots, K \}} \mL^\text{train}_{\theta_k}
\ee
and final best approximate $u_{\theta_K^+}$. Still, $u_{\theta_k}$ depends on the DNN setting, and the losses 
\be\label{eq:losses}\rev{\text{loss}^\text{train} := \mL_{\theta_K^+}^\text{train}} \quad \text{and} \quad \rev{\text{loss} \equiv \text{loss}^\text{test} : =}\mL_{\theta_K^+}^\text{test}
\ee 
can be considerable, even for large $K$. As a consequence, we tune the PINNs in order to find optimal configurations.

\section{Bayesian HPO}\label{sec:HPO}Following \cite[Section 4.3]{chen2021learning}, we employ Bayesian HPO to seek for the optimal PINNs configuration. In this work, we opt for optimizing five HPs namely:
\begin{enumerate}
	\item Learning rate $\alpha$;
	\item Width $N$: number of nodes per layer;
	\item Depth $L-1$: number of dense layers;
	\item Activation function $\sigma$;
	\item Weights for the boundary term error $w_\Gamma$.
\end{enumerate}
Notice that we restrict to constant-width DNNs, i.e.~$N_l=N$ for $N=1, \cdots, L-1$. The search space is $\Lambda$, the Cartesian product of all HP ranges. Every $\lambda \in \Lambda$ writes as:
$$
\lambda = [\alpha, N, L-1, \sigma, w_\Gamma].
$$
The HPO of PINNs can be represented by a bi-level optimization problem \cite{franceschi2017bridge}:
\be
\lambda^\star = \argmin_{\lambda \in \Lambda} \mL^\text{test}_{\theta^\star}[\lambda]\quad \text{with}\quad \theta^\star = \min_{\theta \in \Theta} \mL^\text{train}_\theta [\lambda].
\ee
%or equivalently, as finding:
%$$
%\min_{ \lambda \in \Lambda}\min_{\theta \in \Theta} \mL_\theta[\lambda].
%$$
In the same fashion as with the inner loop in \Cref{picture:PINNs}, we aim at applying an optimizer to the outer loop. Beforehand, notice that:
\begin{itemize}\item $\dim(\Lambda)=5$. Grid search becomes impractical, as a $10$ points per dimension grid would lead to $10^5$ outer loops;
\item The inner loop in \Cref{picture:PINNs} is potentially computationally intensive and time consuming, justifying a more subtle approach than (pure) random search.
\end{itemize}
A natural choice is to apply GP-based Bayesian HPO to obtain approximates to $\lambda^\star$. This method is a good trade-off between brute force grid search and random search. \rev{It consists in applying a GP-based regressor to the previous configurations in order to predict the best next HP setting.} In \Cref{alg:HPO} we summarize the HPO. 
\begin{algorithm}[h!t]
\caption{GP-based Bayesian HPO \cite{yu2020hyper} for $\lambda$}
\begin{algorithmic}
\label{alg:HPO}
\REQUIRE $g_\text{acq}$ (acquisition function)\rev{, $\text{loss}$ (loss function)}
\REQUIRE $\lambda_0$ (initial HPs)
\FOR{$m=0$ \TO $M-1$} 
\STATE Apply GP regression to $(\lambda_i, \rev{\text{loss}}(\lambda_i))$, $i=0,\cdots, m:$\\
$\Rightarrow \rev{\text{loss}}(\lambda)\sim \text{GP}(\mu(\lambda), k(\lambda,\lambda'))$
\STATE $\lambda_{m+1} \leftarrow \argmin_\lambda g_\text{acq}(\mu (\lambda), k(\lambda,\lambda'))$
\STATE Sample $\rev{\text{loss}}(\lambda_{m+1})$
\ENDFOR
\ENSURE $(\lambda_M, \rev{\text{loss}}(\lambda_M))$
\end{algorithmic}
\end{algorithm}
The acquisition function in \Cref{alg:HPO} allows to define the next point at each iteration $m$. \rev{For the sake of clarity, $m$ is coined as an iteration, and $k$ as an epoch, respectively.} A common acquisition function is the negative expected improvement, namely:
\be\label{eq:EI}
- EI(\lambda) = - \IE [\max (f(\lambda)- f(\lambda_m^+),0 )]
\ee
where $\lambda_m^+$ is the best point observed so far and $\IE$ is the expected value. \rev{Concerning the GP, one generally resorts to constant $\mu(\lambda)$ and to Matérn or squared exponential covariance kernel $k(\lambda,\lambda')$\footnote{\url{https://scikit-learn.org/stable/modules/generated/sklearn.gaussian_process.kernels.Matern.html}}. Categorical data such as the activation function $\sigma$ can be transformed into continuous variables (see for example the setting later on in \Cref{subsec:Methodology}).}

The outer loop is conducted throughout an $M$-step GP-based Bayesian HPO. For each $\lambda \in \Lambda$, we define the loss function as being the best (inner) ADAM epoch for $k \in [0,K]$ \rev{with $\text{loss}$ as in \eqref{eq:losses}}:
\be\label{eq:Accuracy}
\rev{\text{loss}}[\lambda] =  \mL_{\theta}^\text{test}[\lambda] \quad \text{for} \quad \theta = \argmin_{k = 0,\cdots, K} \mL^\text{train}_{\theta_k} [\lambda].
%l(\lambda):= \mL_{\theta}^\text{test}[\lambda]\quad \text{for} \quad \theta = \argmin_{k = 0,\cdots, K} \mL^\text{train}_{\theta_k} [\lambda].
\ee
The bi-level optimization produces the final solution $u_{\theta_K^+}[\lambda_M^+]$. We are now ready to perform the numerical experiments.

\section{Numerical experiments}\label{sec:NumExp} In what follows, we set up the numerical experiment's protocol. 
\subsection{Methodology}\label{subsec:Methodology}Throughout, we tune PINNs via HPO defined previously in \Cref{sec:HPO}. Our focus is the Helmholtz equation in two and three dimensions. This has been implemented in the open-source PINNs library DeepXDE 1.1.14 \cite{lu2021deepxde}. Simulations were performed in single float precision on a AMAX DL-E48A AMD Rome EPYC server with 8 Quadro RTX 8000 Nvidia GPUs---each one with a 48 GB memory. DeepXDE supports TensorFlow 1.x, TensorFlow 2.x and PyTorch as backends, and is with single GPU acceleration. \rev{GP-based Bayesian HPO} is performed with the novel high-level HPO library HPOMax\footnote{\url{https://github.com/pescap/HPOMax/}} built upon Scikit-Optimize\footnote{\url{https://github.com/scikit-optimize/scikit-optimize/}}. For reproducibility purposes, results are performed with global random seed $11$.

In this section, we solve the Helmholtz equation in $D= [0,1]^d$, $d=2,3$. For any pulsation $\omega \in \IN$, $\omega \geq 1$, and $\kappa: = 2\pi \omega$, the problem reads:
\begin{align*}
-\Delta u - \kappa^2 u = f \quad \text{in}\quad D,\\
 \mB[u] = 0 \quad \text{on}\quad \Gamma.
\end{align*}
wherein $d,f$ and $u$ are summarized for both the Dirichlet and Neumann cases in \Cref{tab:OverviewNumericalExperiments}. For Dirichlet case, we enforce hard constraint BCs by applying the transformation:
$$
\hat{u}_\theta(\bx) = (1-x^2)(1-y^2)u_\theta(\bx)
$$
to the net. \rev{Notice that the usage of this transformation is common, but that other choices are possible as long as they enforce the BCs. Finally, a comprehensive tutorial was added to the DeepXDE documentation\footnote{\url{https://deepxde.readthedocs.io/en/latest/demos/pinn_forward/helmholtz.2d.dirichlet.hpo.html}}}.

\begin{table}[H]
\renewcommand\arraystretch{1.5}
\begin{center}
\footnotesize
%\resizebox{8.5cm}{!} {
\begin{tabular}{
    |>{\centering\arraybackslash}m{2cm}% instead of "p" is "m"
    |>{\centering\arraybackslash}m{.5cm}
    |>{\centering\arraybackslash}m{2.5cm}
    |>{\centering\arraybackslash}m{2.5cm}
    |>{\centering\arraybackslash}m{.5cm}|
    }
    \hline
Case & d & $f$ & $u$\\ \hline\hline
Dirichlet & 2 &$\kappa^2 \sin(\kappa x) \sin (\kappa y)$ &$\sin(\kappa x) \sin(\kappa y)$ \\\hline 
Neumann & 3 & $2 \kappa^2 \cos(\kappa x) \cos(\kappa y)$ &$\cos(\kappa x) \cos(\kappa y)$ \\\hline 
 \end{tabular}%}
\end{center}
\caption{Overview of the data for the two cases considered throughout this section.}
 \label{tab:OverviewNumericalExperiments} 
\end{table}

The training (resp.~testing) collocation points are generated randomly with a precision of $r=10$ (resp.~$r=30$) points per wavelength per dimension, with $\frac{2\pi}{\kappa}$ the wavelength. The inner loop is solved via ADAM over $K=50{,}000$ epochs with Glorot uniform initialization \cite[Chapter 8]{bengio2017deep}. The HPO is run over $M=100$ iterations, with negative expected improvement in \eqref{eq:EI} as the acquisition function. The default parameters are set to:
\begin{align*}
\text{Dirichlet:}  \quad \lambda_0& = [10^{-3}, 4, 50, \sin],\\
\text{Neumann:}   \quad \lambda_0 &= [10^{-3}, 3, 275, \sin , 400].
\end{align*}
In \Cref{tab:OverviewHP}, we summarize notations for HPs. We list the associated symbol and range for each HP. \rev{We use the default parameters of the gp$\_$minimize Scikit-optimize function (more details in documentation\footnote{\url{https://scikit-optimize.github.io/stable/modules/generated/skopt.gp_minimize.html}}). Regarding the search space, integer variables are processed as continuous variables in the HPO. Furthermore, 
 categorical variables are converted into normalized binary vectors, and applied an argmin function upon. For example, $\sigma \in \{\tanh,\sin\}$ is transformed into a two-dimensional vector with $[1,0]$ as $\sigma=\tanh$ and $[0,1]$ as $\sigma=\sin$. Therefore, if the surrogate prediction is $[0.4,0]$, one will use a $\tanh$ function at next iteration. Additionally, the first $10$ iterations are performed with random points generation.}
\begin{table}[h!t]
\renewcommand\arraystretch{1.5}
\begin{center}
\footnotesize
%\resizebox{11.5cm}{!} {
\begin{tabular}{
    |>{\centering\arraybackslash}m{3.5cm}% instead of "p" is "m"
    |>{\centering\arraybackslash}m{1cm}
    |>{\centering\arraybackslash}m{2.6cm}
    |>{\centering\arraybackslash}m{2.6cm}
    |>{\centering\arraybackslash}m{1.4cm}|
    }
    \hline
\multirow{2}{*}{Hyper-parameter} & \multirow{2}{*}{Symbol} & \multicolumn{2}{|c|}{Range} & Log-\\ \cline{3-4}
&&Dirichlet &  Neumann& transform\\ \hline \hline
Learning rate& $\alpha$ & [$10^{-4}, 5\times 10^{-2}$] & [$10^{-5}, 5\times 10^{-2}$] & Yes\\\hline 
Width  (\# dense nodes)& $ N$ & $[5,500]$ & $[5,500]$ & No \\ \hline 
Depth (\# hidden layers)& $L-1$ & $[1,10]$ & $[1, 5]$ & No \\\hline 
Activation function & $\sigma$ &  \multicolumn{2}{c|}{$\{\sin, \sigmoid, \tanh\}$} &-\\ \hline
Weight  & $w_\Gamma$ &  - & $[1, 10^{7}]$ & Yes \\\hline
 \end{tabular}%}
\end{center}
\caption{Overview of the \rev{HPs}. We provide the default values for HPO. Dirichlet BCs in \eqref{eq:Hard1}-\eqref{eq:Hard2} remove the boundary loss term, and the weight $w_\Gamma$.}
\label{tab:OverviewHP} 
\end{table}

\subsection{Dirichlet case: HPO}\label{subsec:dirichletHPO}
To begin with, we focus on the Dirichlet case. We set $d=2$, $\omega=2$---i.e~$\kappa = 4\pi$. The collocation points are with:
$$
|\mT^\text{train}| =400 \quad \text{and}\quad |\mT^\text{test}|  = 3{,}600. 
$$
HPO leads to an optimal configuration at step $97$ where
\be\label{eq:BestDirichlet}
\lambda_M^+ = \lambda_{97} = [10^{-4}, 2, 275, \sin]
\ee
with
\be 
\rev{\text{loss}}[\lambda_{97}] = 1.13 \times 10^{-3} \quad \text{and}\quad |\Theta_{97}| = 77{,}001 .
\ee
\begin{figure}[!htb]
\center
\includegraphics[width=.85\textwidth]{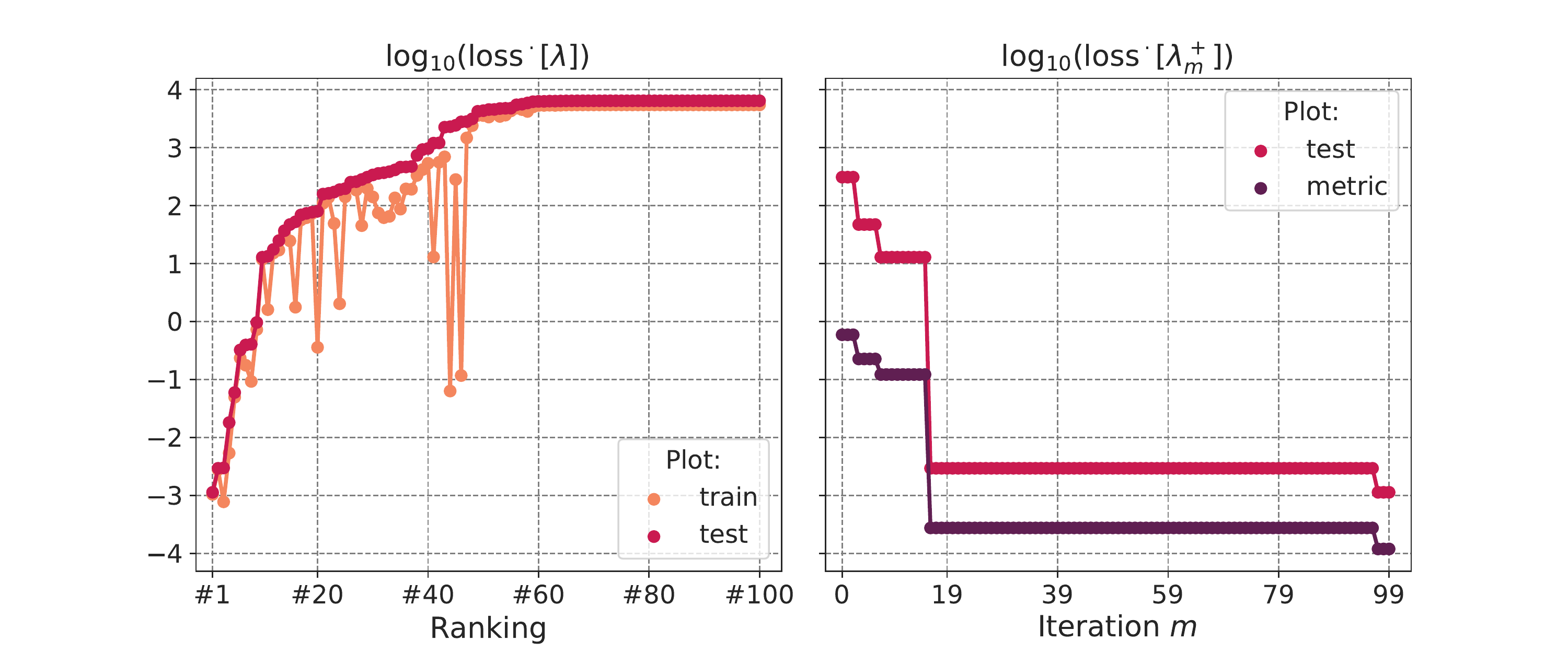}
\caption{\rev{HPO: Ordered $\text{loss}$ in ascending order (left) and best $\text{loss}$ at iteration $m$ (right).  Remark that the loss ranges between $10^{-3}$ and $10^{4}$. Results are in log scale.}}
\label{fig:HPO1}       
\end{figure}
\rev{On the left side of \Cref{fig:HPO1}, we show $\{\text{loss}^\cdot[\lambda_m]\}_{m=0}^{99}$ for $\cdot \in \{\text{train},\text{test}\}$ ordered in ascending order. It allows to analyze the distribution of $\text{loss}^\cdot$, $\cdot \in \{\text{train},\text{test}\}$, over the iterations. On the right side \rev{of \Cref{fig:HPO1}}, we represent the best $\text{loss}[\lambda]$ at step $m$. Acknowledge that the latter consists in the output of the HPO at step $m$ (refer to \Cref{alg:HPO}.) We observe that the training and testing errors are close together for the best and worst configurations, although some of them do not accomplish this pattern in the middle part of the figure (for instance, iteration $11$ is with \rev{$\text{loss}^\text{train}[\lambda_{11}] =1.60$} and $\rev{\text{loss}^\text{test}[\lambda_{11}]} = 13.2$).} Furthermore, only few $\lambda_m$ lead to small loss. To illustrate this, it is remarked that the 80\% worst configurations are with more than $10^2$ as test error. Also, the loss ranges from $10^{-3}$ to $10^{4}$, which implies a high variability in the training results. These two remarks justify the core of the tuning procedure. Next, we analyze the performance of HPO along iteration $m$ for $m=0,\cdots,M-1$. On the right side, we show HPO best iterate at step $m$. First, notice the correlation between the test and metric curves. In this simulation, an improvement in the test loss comes with improvement in the solution accuracy. Second, the test loss is with a jump of more than $3$ orders of magnitude at step $16$ ($\rev{\text{loss}[\lambda_{16}]}= 2.94\times 10^{-3}$ and $\rev{\text{loss}^\text{metrics}}[\lambda_{16}]= 2.75\times 10^{-4}$). Next, the best configuration is found at step $97$ ($\rev{\text{loss}}[\lambda_{97}]= 1.13\times 10^{-3}$ and $\rev{\text{loss}^\text{metric}}[\lambda_{97}]= 1.19\times 10^{-4}$). To reformulate, a good configuration was found at step $16$, and a slightly more optimal (by around a factor of $2$) emerges at step $97$. Concerning execution times, HPO takes $10.0$ hours. Each iteration $m$ takes $6$ minutes on average. The standard deviation is of $3$ minutes and $8$ seconds. The number of trainable parameters is $4.23 \times 10^5 \pm 6.38 \times 10^5$ (average $\pm$ standard deviation).
\begin{figure}[!htb]
\center
\includegraphics[width=.6\textwidth]{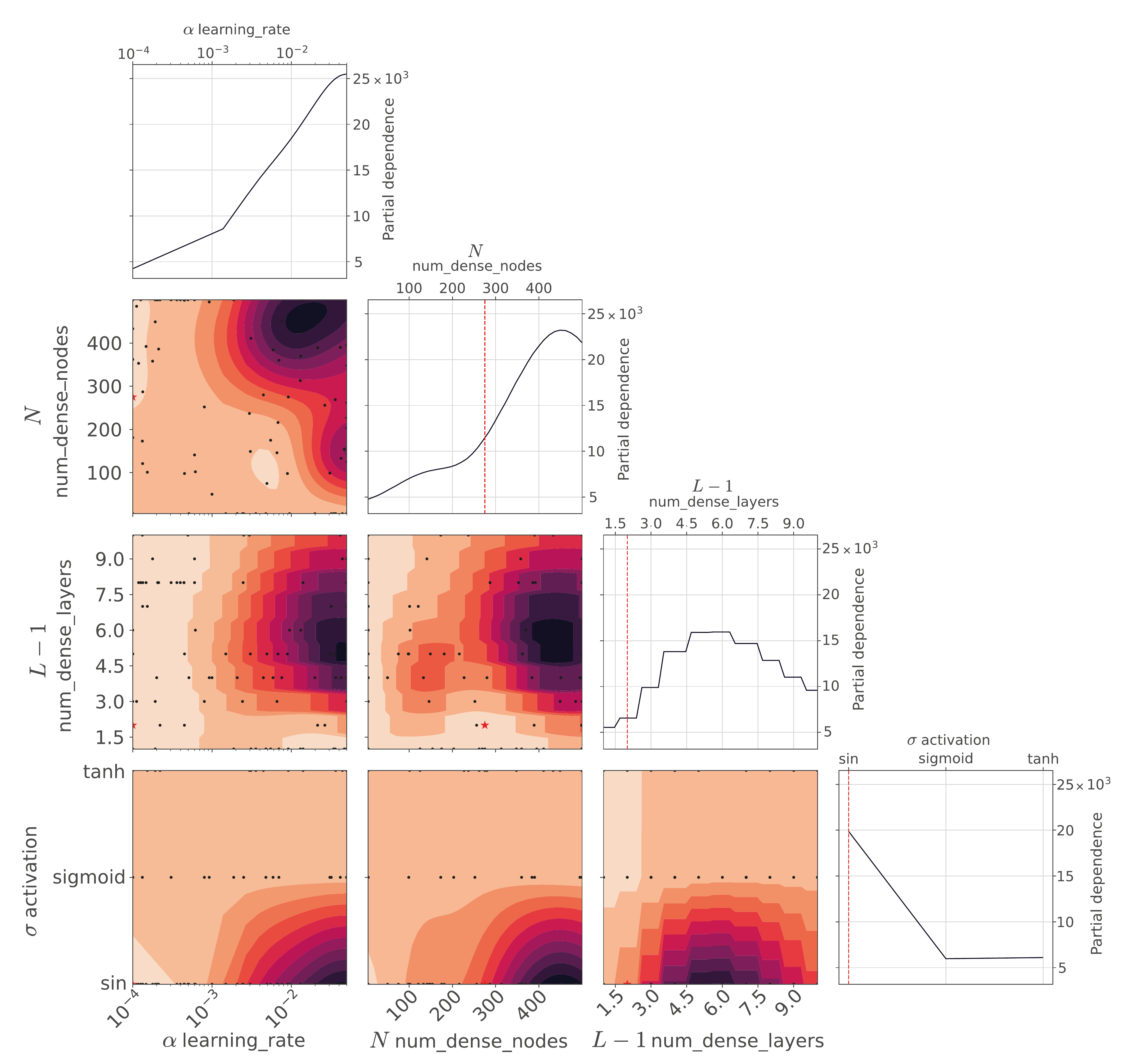}
\caption{Partial dependence plots of the objective function.}
\label{fig:HPO3}
\end{figure}

%Rpartial dependence plot illustrates the estimated relative importance of each variable on the loss; if the partial dependence is constant, it indicates that the variable has no effect on the loss

Next, we focus on the HPs. In \Cref{fig:HPO3}, we represent the partial dependence plots of the objective function. \rev{It estimates the relative dependence of each dimension on the loss (after averaging out all other dimensions).} We remark that the partial dependence is bounded away from a constant for all HPs, testifying for a complex loss function depending on all the HPs. We remark that some areas (in light salmon color) show better behaviors. To sum up, the plot suggests using a small learning rate. For the number of dense nodes, several zones appear, but we focus on the area for the minimum. It suggests to use:
\be\label{eq:suggestion} 
\alpha\approx 10^{-4}, \quad N \in [250,500],\quad L-1  \leq 3, \quad \sigma=  \sin.
\ee
Notice in the last row that $\sigma =\sin$ leads to the best and the worst error losses \rev{with respect to $\alpha$ and $N$.}
\begin{figure}[!htb]
\center
\includegraphics[width=.35\textwidth]{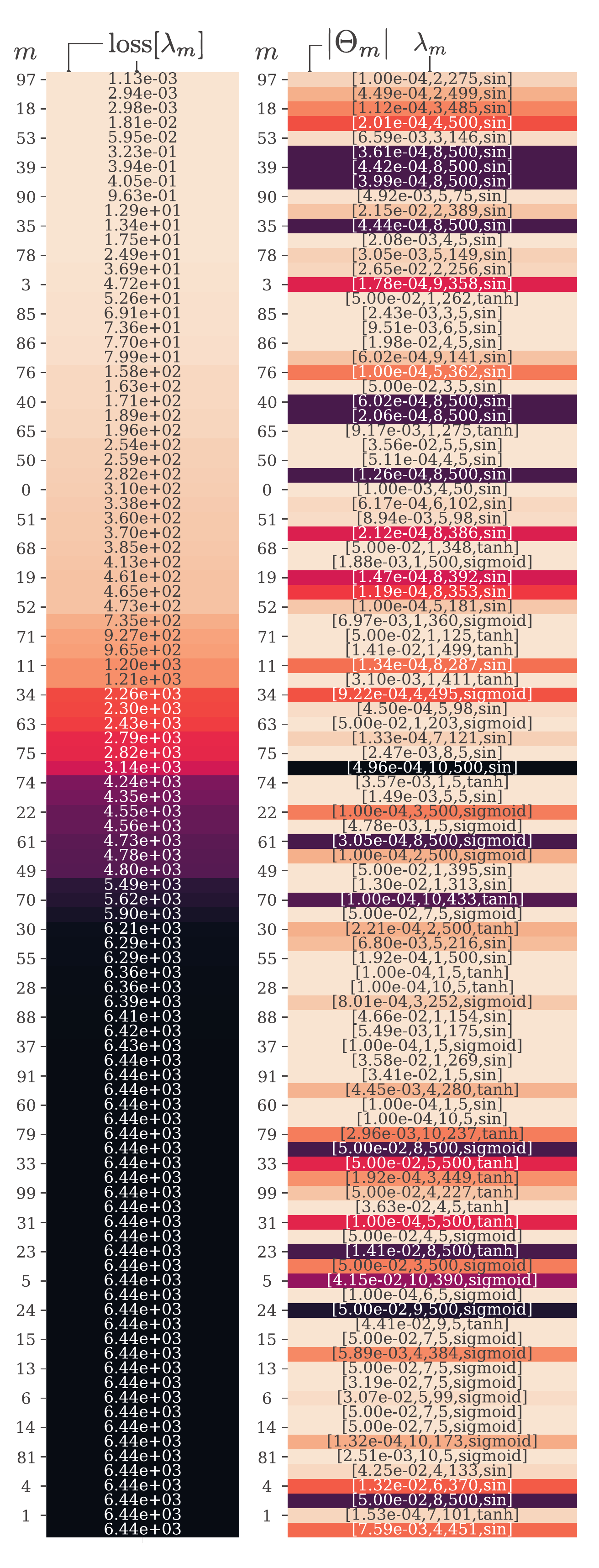}
\caption{Summary of HPO and error. Left column shows the ranking for $\text{loss}[\lambda]$. Right column portrays the HPs, and the color is for the magnitude of $\Theta$.}
\label{fig:HPO4}       
\end{figure}

To finish the analysis, we scrutinize the top configurations more into detail. In \Cref{fig:HPO4}, we show a global overview of HPO in a manner that allows to observe the HPOs, the loss and number of trainable parameters. Again, we remark that very few configurations are with an error lesser than $1$. Now, we can find patterns in the best configurations. To begin with, notice that the $10$ best configurations are with $\sin$ activation function. Furthermore, some configurations are with heavy computational cost, due to the DNN with e.g.~$10 \times 500$ or $8 \times 500$ deep$\times$width. As inferred in \Cref{fig:HPO4}, the best configurations are with learning rates around $10^{-4}$ and few layers $L-1$. As hinted before, several configurations do not converge at all, justifying HPO.

\subsection{Dirichlet case: Best configuration}\label{subsec:DirichletBest}
We focus on the best iterate $\lambda_M^+ = \rev{\lambda_{97}}=[10^{-4}, 2, 275, \sin]$.
\begin{figure}[!htb]
\center
\includegraphics[width=.45\textwidth]{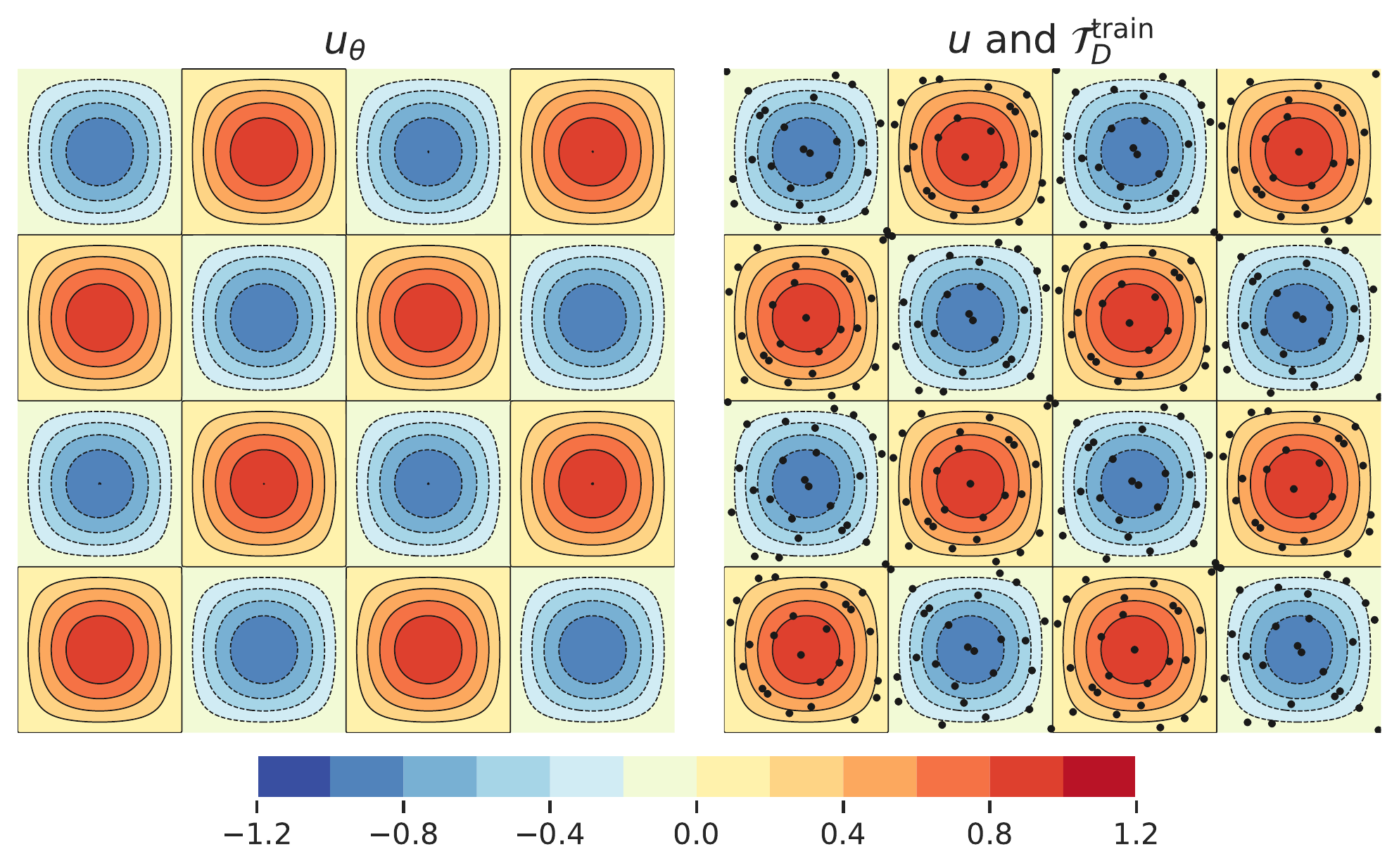}
\caption{Solution obtained by PINNs $u_\theta$ (left) v/s exact solution $u$ and the collocation points $\mT^\text{train}_D$ (right).}
\label{fig:B10}       
\end{figure}
In \Cref{fig:B10}, we represent $u_\theta$ the solution obtained by the PINN (top). We compare it to the exact solution $u$ (bottom). We also represent the $400$ training collocation points $\mT^\text{train}_D$, which allow to define the loss function. Remark that $u_\theta$ and $u$ coincide.

\begin{figure}[!htb]
\center
\includegraphics[width=.45\textwidth]{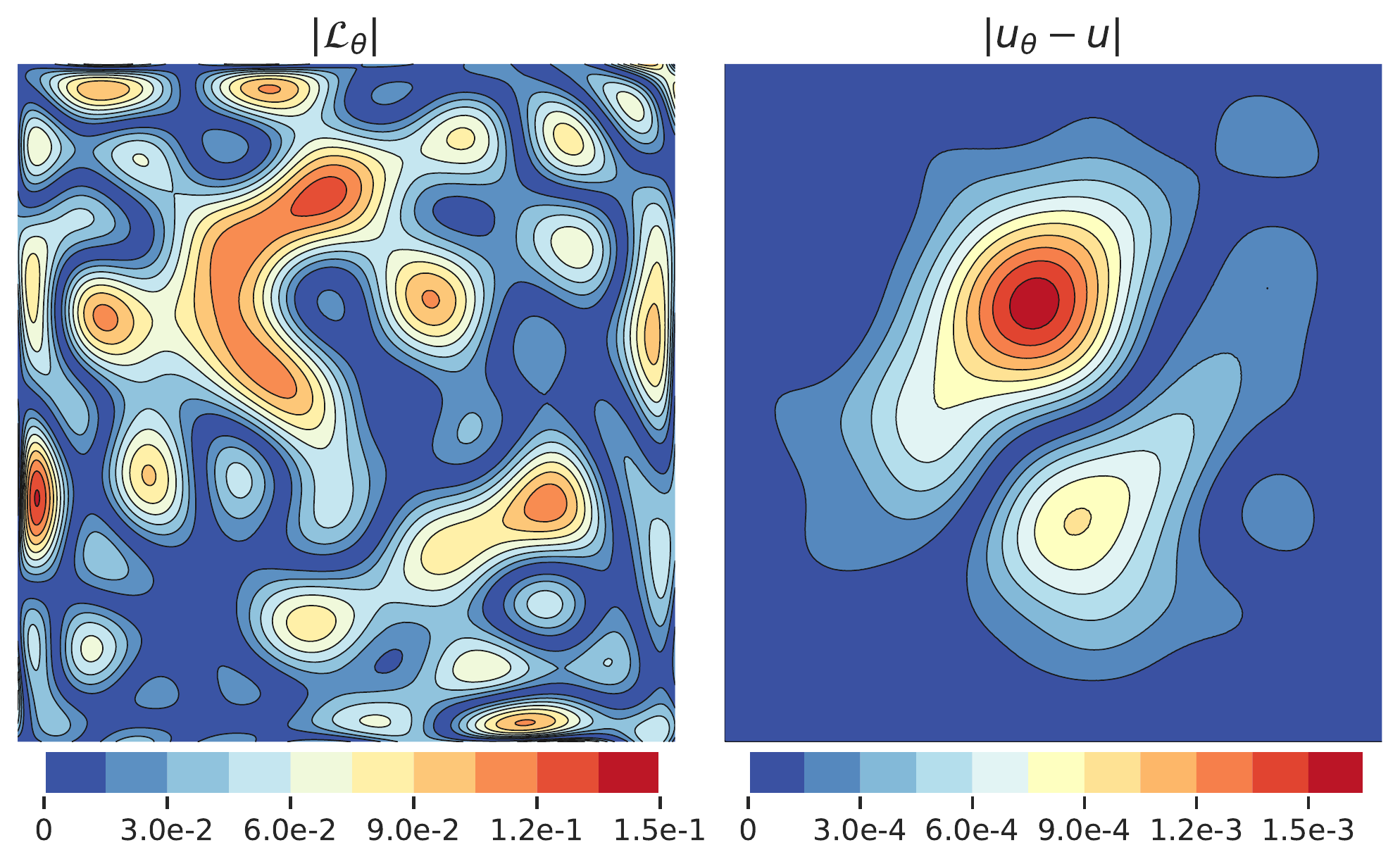}
\caption{Pointwise loss $|\mL_\theta|$ and residual error $|u_\theta-  u|$.}
\label{fig:B11}       
\end{figure}

\begin{figure}[!htb]
\center
\includegraphics[width=.45\textwidth]{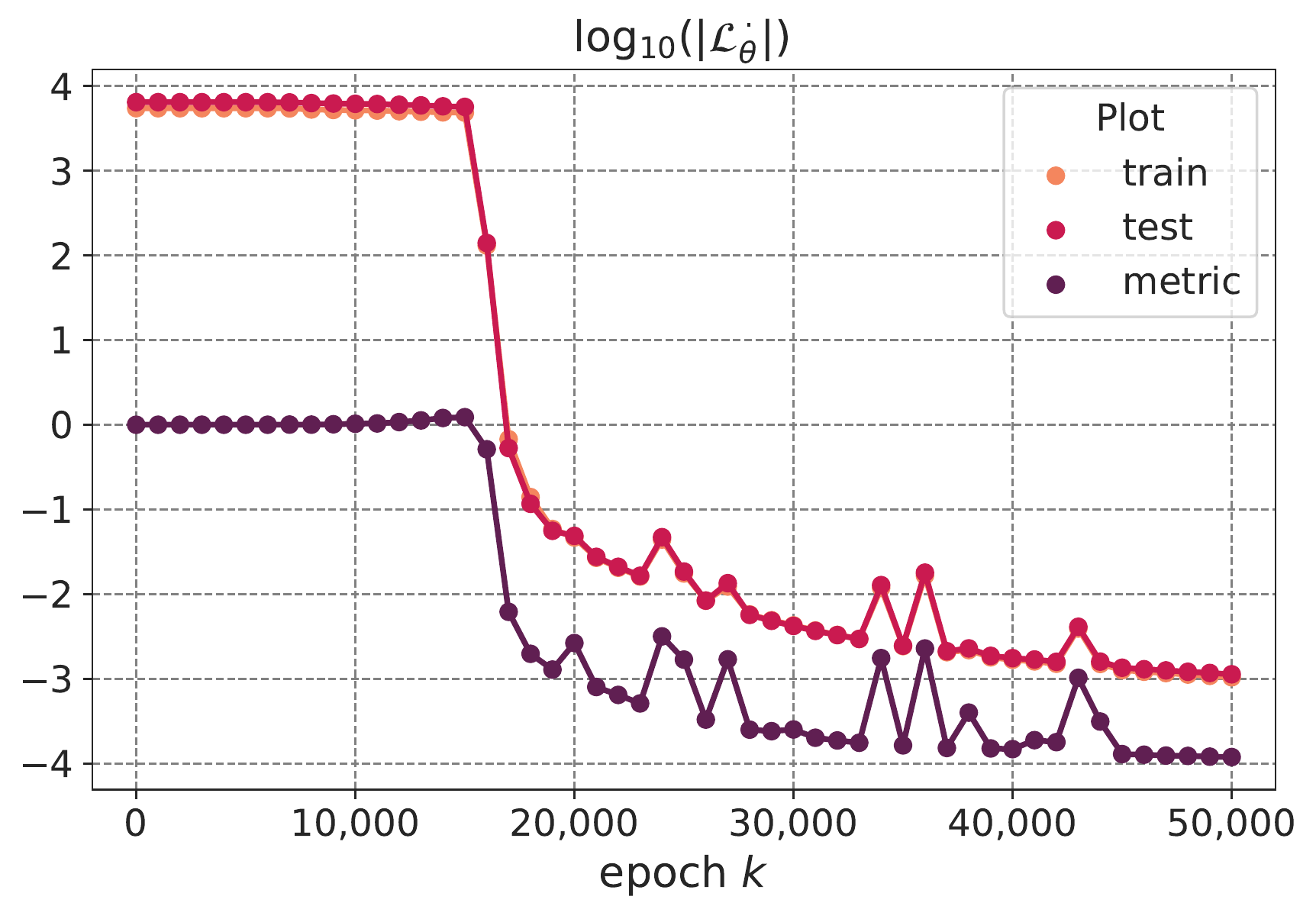}
\caption{Convergence results for the ADAM optimizer for $\lambda_{97}$. Best iterate is for $k=K=50{,}000$.}
\label{fig:B12}       
\end{figure}

For the sake of completeness, we showcase the pointwise loss in \Cref{fig:B11} (left).  Furthermore, we also represent the pointwise residual error $|u_\theta - u|$ (right). The pointwise error ranges between $0$ and $1.6 \times 10^{-3}$. Remark that the largest error for $|u_\theta-u|$ concentrates away from the boundary.
Furthermore, we represent the convergence for the ADAM optimizer in \Cref{fig:B12}. The optimal test error---as well as for testing and metric errors---is obtained at epoch $k=50{,}000$. Remark that the training and testing error coincide for all epochs. Also, the metric is highly correlated with the testing error.
The execution time is $196.6$ seconds. 

\subsection{Dirichlet case: $(h,\kappa)$-analysis}\label{subsec:hk}
\begin{figure}[!ht]
\center
\includegraphics[width=.6\textwidth]{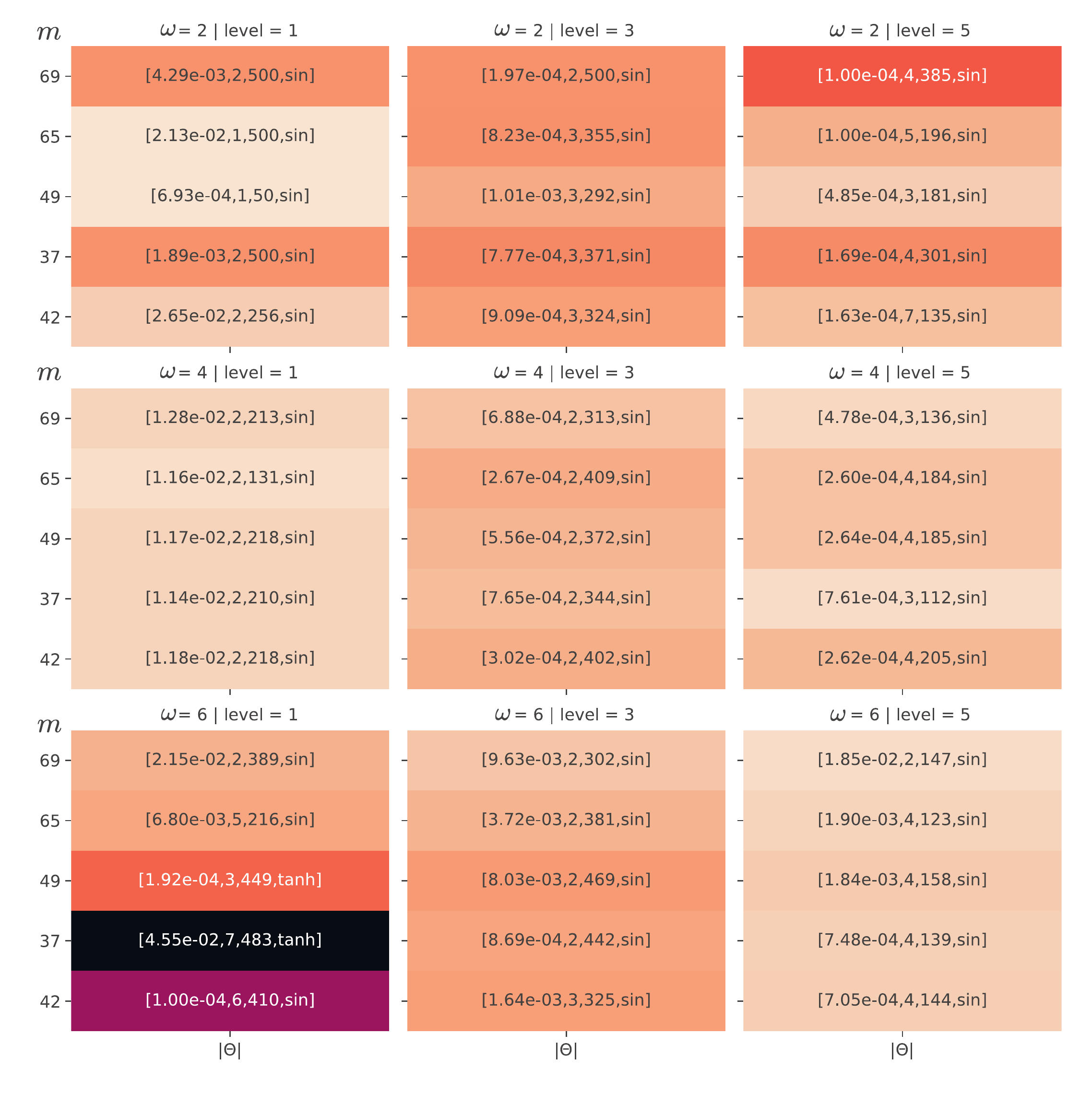}
\caption{Summary of the Top 5 configurations for each level (row) and pulsation (column). The color represents the number of trainable parameters $|\Theta|$. Notice that the observations in \eqref{eq:suggestion} remain valid.}
\label{fig:hk2}       
\end{figure}

We complete the analysis for the Dirichlet case by performing an $(h,\kappa)$-analysis. To that extent, we solve HPO with the same setting as in \Cref{subsec:dirichletHPO}, for $\omega = \{2,4,6\}$ and levels $\{\mT_1,\mT_3,\mT_5\}$ with
$$
|\mT_l| = n_l^2\quad \text{with} \quad n_l := 10 \times 2^{l - 1} \quad \text{for}\quad l= \{1,3,5\}.
$$
To sum up, $\mT_1,\mT_3,\mT_5$ correspond to setting $n_x = 10$, $40$, $160$ points per dimension. This leads to
$$
|\mT_1| = 100,\quad |\mT_3| = 1{,}600, \quad |\mT_5| = 25{,}600.
$$
The configurations and their respective precision $r$ are depicted in \Cref{tab:Overviewhk}.
\begin{table}[!htb]
\renewcommand\arraystretch{1.3}
\begin{center}
\footnotesize
%\resizebox{7cm}{!} {
\begin{tabular}{
    |>{\centering\arraybackslash}m{1cm}% instead of "p" is "m"
    |>{\centering\arraybackslash}m{1.5cm}
    |>{\centering\arraybackslash}m{1.5cm}
    |>{\centering\arraybackslash}m{1.5cm}|
    }
    \hline
$\omega $& $\mT_1$ & $\mT_3$ & $\mT_5$\\ \hline\hline
$2$ & $ 5.0 $ & $20.0$ & $80.0$\\\hline 
$4$ & $2.5$ & $10.0$ & $40.0$\\ \hline 
$6$ & $1.7$ & $6.7$ & $26.7$ \\\hline 
 \end{tabular}%}
\end{center}
\caption{Precision $r$ for all the configurations. For example, the top-left cell means that for pulsation $\omega =2$ and level $\mT_1$, one has a precision of $r=5.0$ points per wavelength per dimension. Notice that a rule of thumbs in FEM community is using $r=10$.}
\label{tab:Overviewhk} 
\end{table}
\begin{figure*}[!htb]
\center
\includegraphics[width=.9\textwidth]{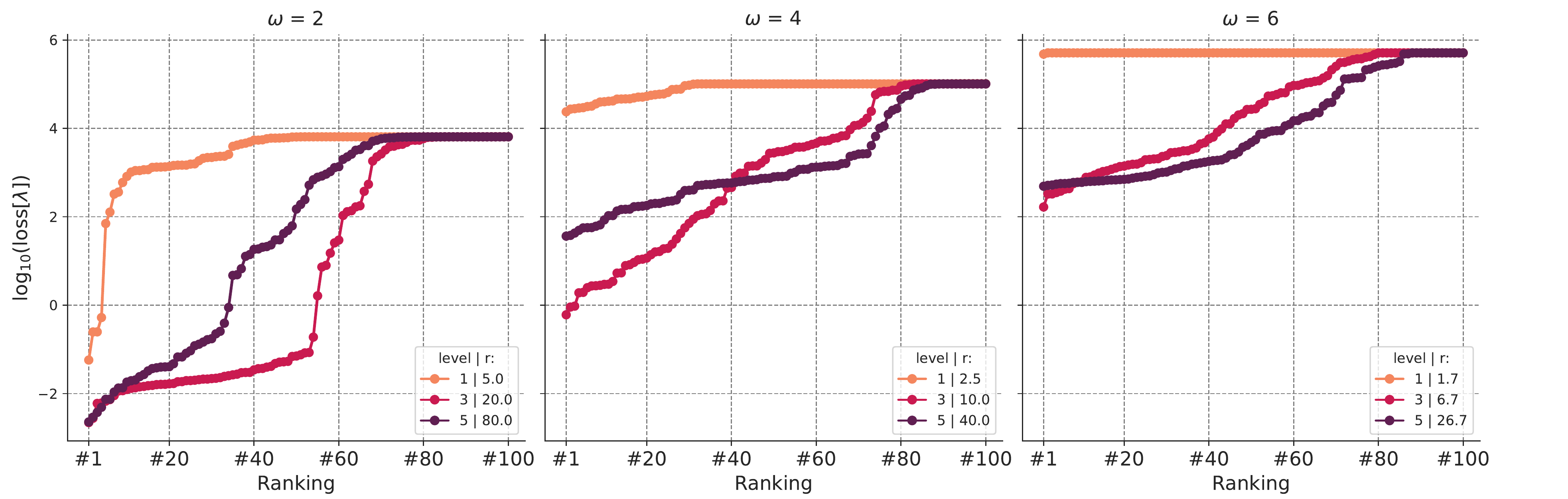}
\caption{Ordered values of $\text{loss}[\lambda]$ function to the levels (1, 3, 5) and the pulsation $\omega =2$ (left), $\omega =4$ (center) and $\omega= 6$ (right).}
\label{fig:hk0}       
\end{figure*}
In \Cref{fig:hk0}, we plot the result of HPO for all those configurations.
\begin{figure*}[!htb]
\center
\includegraphics[width=.9\textwidth]{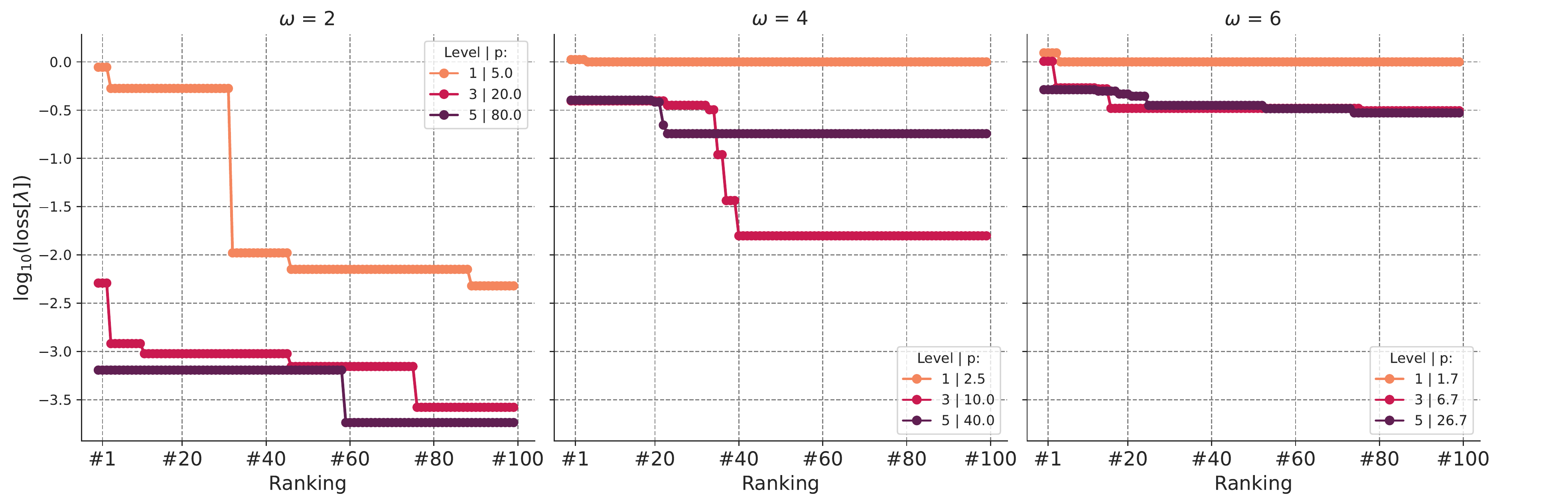}
\caption{\rev{$\text{loss}^\text{metric}[\lambda_m]$} function to the levels (1, 3, 5) and the pulsation $\omega =2$ (left), $\omega =4$ (center) and $\omega= 6$ (right).}
\label{fig:hk1}       
\end{figure*}
For $\omega=2$, we remark that increasing level leads to better results. Also, all levels lead to satisfactory ranges. A important statement is that results deteriorate with increasing frequencies (from left to right). For $\omega=4$ and $\omega=6$, the best results are obtained with $\mT_3$, while $\mT_1$ does not converge (due to the small $r$). In the same fashion, we represent the metric $\rev{\text{loss}}^\text{metric}[\lambda_m]$ along HPO in \Cref{fig:hk1}. Here, we observe that for $\omega=4$, the optimal metric is obtained for $\mT_3$, despite not showing the best loss in \Cref{fig:hk0}. The deterioration with $\omega$ is made clear here. We see that the results for $\omega=6$ are not satisfactory, despite showing an improvement during HPO. For the sake of completeness, we sum up the Top 5 configurations for all the settings in \Cref{fig:hk2}. Notice that the HPs recommendation in \eqref{eq:suggestion} remain valid. 

\rev{It is surprising to notice that higher levels do not necessarily lead to more successful HPO. Acknowledge that for a given $\lambda$, some results exist concerning the generalization error according to $|\mT^\text{train}|$ and the training error \cite[Theorem 2.6]{Mishra2020EstimatesOT}. Yet, in our setting we are varying both $\lambda$ and $|\mT^\text{train}|$ at the same time. Our results hint at understanding how these parameters interact more into detail, but are out of the scope of this manuscript.}

\subsection{Dirichlet case: Comparison to FEM}To finish, we compare the best configuration for each setting to FEM. The FEM solutions are obtained with DOLFINx \cite{logg2010dolfin} with 10 points per wavelength, and piecewise polynomial elements of order $2$. We use both a direct solver (LU factorization) and iterative solver (unpreconditioned GMRES with relative tolerance $10^{-5}$, no restart and maximum number of iterations $10{,}000$). The results are summarized in \Cref{tab:FinalFEM}. We start with accuracy (second column ``metric''). We remark that tuned PINNs outperform FEM in terms of accuracy for $\omega = 2$ and for $\omega=4,\mT_3$ (in bold notation). Remark that the error decreases with the level for FEM, by virtue of the quasi-optimality result for FEM for high enough precision \cite{ern2004theory}. Still, the tuned PINNs follow this pattern only for small frequency $\omega =2$. Then, the best accuracy is obtained on $\mT_3$ for $\omega \in \{ 4,6\}$. Next, we focus on the third column with execution times. Both FEM alternatives outperform PINNs in all cases, by around a factor of $100$. In this case, the fastest method is LU factorization. These results show the limitation of PINNs with the frequency, and pave the way toward working on frequency stable schemes.
\begin{table}[!htb]
\renewcommand\arraystretch{1.3}
\begin{center}
\footnotesize
%\resizebox{8.75cm}{!} {
\begin{tabular}{
    |>{\centering\arraybackslash}m{.5cm}% instead of "p" is "m"
    |>{\centering\arraybackslash}m{1cm}
    |>{\centering\arraybackslash}m{1.4cm}
    |>{\centering\arraybackslash}m{1.4cm}
    |>{\centering\arraybackslash}m{1.4cm}
     ||>{\centering\arraybackslash}m{1.4cm}
    |>{\centering\arraybackslash}m{1.4cm}
    |>{\centering\arraybackslash}m{1.4cm}|
    }
    \hline    
\multicolumn{2}{|c}{Setting}&\multicolumn{3}{|c||}{metric} &\multicolumn{3}{c|}{$t_\text{exec}$ (s)}  \\ \hline
$\omega$&  Case & $\mT_1$ & $\mT_3$ & $\mT_5$ & $\mT_1$ & $\mT_3$ & $\mT_5$  \\ \hline\hline
\multirow{3}{*}{$2$} & PINN  & {\bfseries\boldmath 4.8$\times$10$^{-3}$} & {\bfseries\boldmath 2.6$\times$10$^{-4}$} & {\bfseries\boldmath 1.8$\times$10$^{-4}$ }&
2.0$\times$10$^2$ & 2.4$\times$10$^2$ & 4.7$\times$10$^3$ \\ \cline{2-8} 
& LU& 7.4$\times$10$^{-2}$ &2.2$\times$10$^{-2}$ &5.7$\times$10$^{-3}$ & 7.5$\times$10$^{-3}$ & 5.6$\times$10$^{-2}$ & 2.6$\times$10$^0$ \\ \cline{2-8}
&GMRES& 6.9$\times$10$^{-2}$ & 2.4$\times$10$^{-2}$ & 6.7$\times$10$^{-3}$ & 9.6$\times$10$^{-2}$ & 2.0$\times$10$^{-1}$ & 2.3$\times$10$^1$ \\ \hline \hline
\multirow{3}{*}{$4$} & PINN  & 1.6$\times$10$^0$ & {\bfseries\boldmath 1.6$\times$10$^{-2}$} & 1.8$\times$10$^{-1}$ &
1.8$\times$10$^2$ & 2.0$\times$10$^2$ & 9.8$\times$10$^2$\\ \cline{2-8} 
& LU& 4.0$\times$10$^{-1}$ & 2.0$\times$10$^{-2}$ & 5.6$\times$10$^{-3}$ &  8.6$\times$10$^{-3}$ & 5.4$\times$10$^{-2}$ & 2.6$\times$10$^0$  \\ \cline{2-8}
&GMRES &4.0$\times$10$^{-1}$ & 2.4$\times$10$^{-2}$ & 7.0$\times$10$^{-3}$ & 2.3$\times$10$^{-2}$ & 9.4$\times$10$^{-1}$ & 2.3$\times$10$^1$\\ \hline \hline
\multirow{3}{*}{$6$} & PINN  & 1.2$\times$10$^0$ &3.1$\times$10$^{-1}$ & 3.6$\times$10$^{-1}$ & 1.6$\times$10$^2$ & 2.0$\times$10$^2$ & 6.8$\times$10$^2$\\ \cline{2-8}
& LU& 9.0$\times$10$^{-1}$ & 3.5$\times$10$^{-2}$ & 5.5$\times$10$^{-3}$ & 9.6 $\times$10$^{-3}$ &5.4$\times$10$^{-2}$ & 2.6$\times$10$^0$ \\ \cline{2-8}
&GMRES&  9.0$\times$10$^{-1}$ & 3.5$\times$10$^{-2}$ & 7.0$\times$10$^{-3}$ & 4.4$\times$10$^{-2}$ & 1.1$\times$10$^0$ & 2.3$\times$10$^1$\\ \hline 
 \end{tabular}%}
\end{center}\caption{Training v/s assembly+solver times : PINNs v/s FEM (LU and GMRES). Cases where PINNs outperform FEM are highlighted in bold notation.}
\label{tab:FinalFEM} 
\end{table}

\subsection{Neumann case:HPO}\label{subsec:neumannHPO}We apply the exact same procedure to the $3D$ Neumann case. We define $n_x := 10 r$, $\mT^\text{train}=\mT^\text{test}$, $|\mT_D|= n_x^2$, and
$$
|\mT_\Gamma| := 2^{d-1} d  n_x^{d-1} = 16 n_x^2.
$$
The HPO process leads to:
$$
\lambda_M^+ = \lambda_{92} =  [10^{-4}, 10, 207, \sin,1.0]
$$
along with
$$
\rev{\text{loss}}[\lambda_{92}] = 1.69 \quad \text{and} \quad \mL^\text{metric}_{\theta_K^+}[\lambda_{92}] = 3.7 \times 10^{-1},
$$

\begin{figure}[!htb]
\center
\includegraphics[width=.75\textwidth]{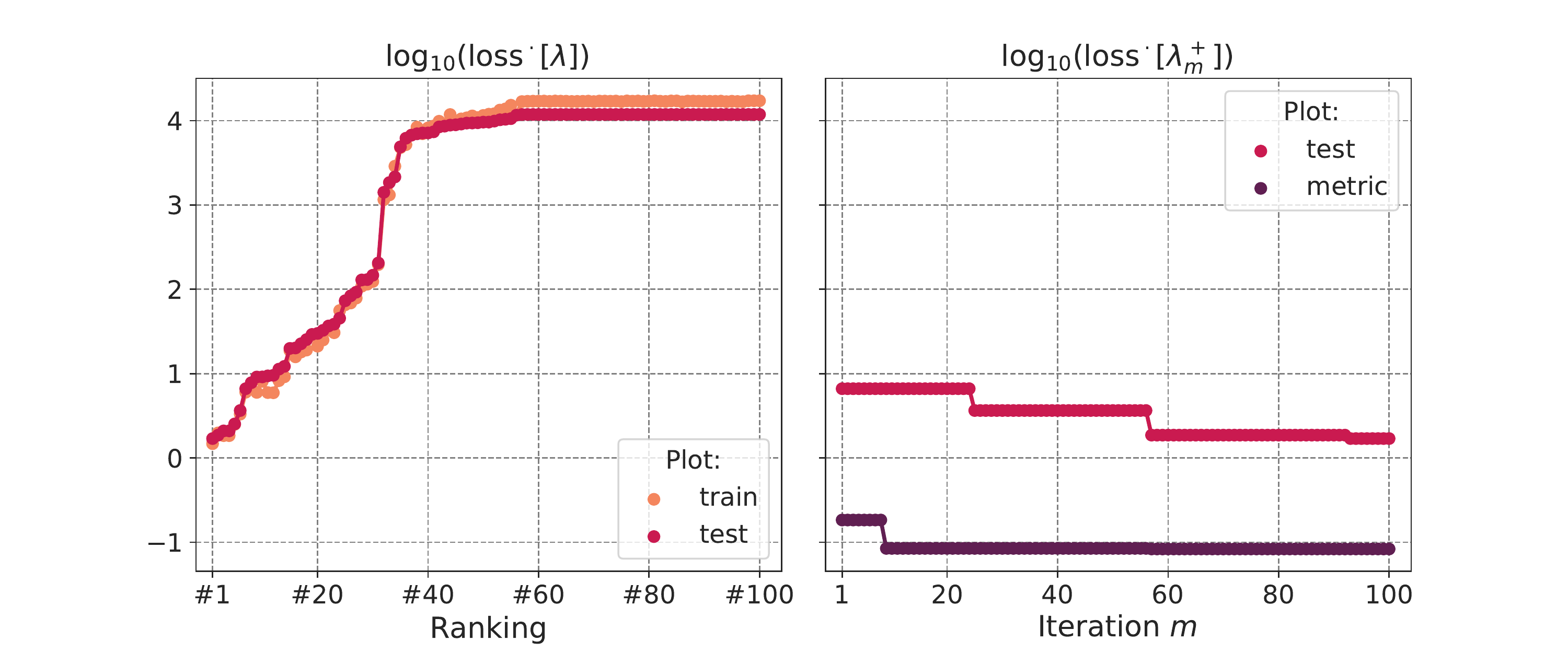}
\caption{HPO: Ordered \rev{$\text{loss}$} in ascending order (left) and best \rev{$\text{loss}$} at iteration $m$ (right). \rev{Results are in log scale.}}
\label{fig:3HPO1}       
\end{figure}

Similar to \Cref{fig:HPO1}, we plot the loss for $m$ in \Cref{fig:3HPO1}. 
Again, the training and test errors are similar. The loss ranges from $1.69$ to $1.19 \times 10^{4}$. Notice that the improvement induced by the HPO (right) is lesser than for the Dirichlet case. This is due to
$$
\lambda_0 = [10^{-3}, 3, 275, \sin, 400]
$$
with
$$
\rev{\text{loss}}[\lambda_0] = 6.62 \quad \text{and} \quad \rev{\text{loss}}^\text{metric}[\lambda_0] = 1.83 \times 10^{-1},
$$
as being an effective configuration. Acknowledge that $\lambda_0$ for the Neumann case was inspired by $\lambda_M^+$ for the Dirichlet case in \eqref{eq:BestDirichlet}. This hints at some similarities between optimal HPs for different problems. Also,  $\rev{\text{loss}}^\text{metric}[\lambda]$ does not reach its minimum value at iteration $m=97$.

\begin{figure}[!htb]
\center
\includegraphics[width=.42\textwidth]{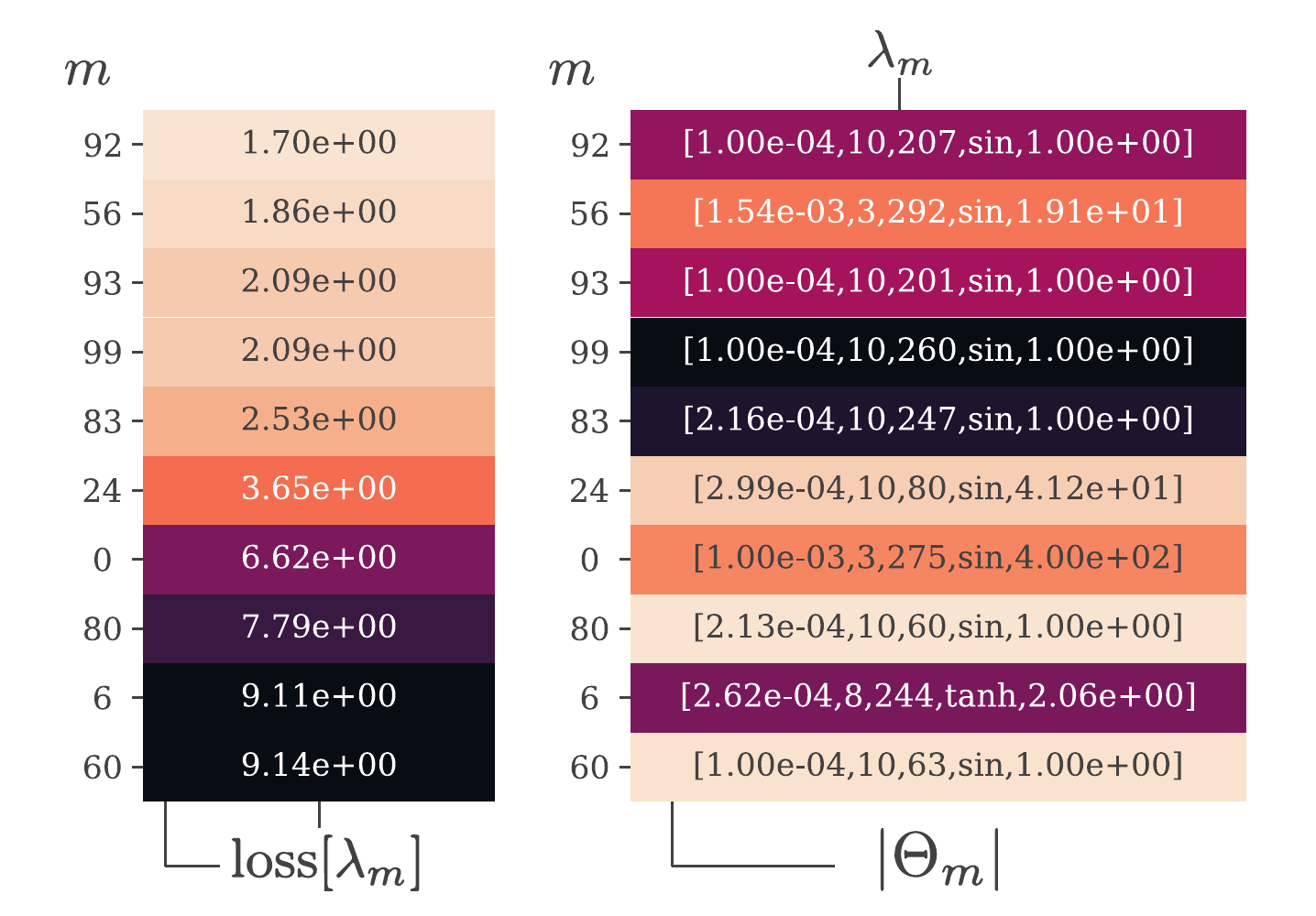}
\caption{Top 10 configurations. Left column shows the ranking for \rev{$\text{loss}[\lambda]$} (color and text). Right column portrays the HPs, and the color is for the magnitude of $\Theta$.}
\label{fig:3HPO4}       
\end{figure}

We scrutinize the top configurations more into detail. In \Cref{fig:3HPO4}, we show a global overview of HPO in a manner that allows to observe the HPs, the loss and number of trainable parameters. 

To finish, we compare PINNs performance to FEM. For PINN, we choose the second configuration:
$$
\lambda_{56} = [1.54 \times 10^{-3}, 3, 292, \sin, 19.1]
$$
as it has fewer trainable parameters. Next, we compare this configuration to FEM. The results are summarized in \Cref{tab:3DFEM}. We define $N_\text{dofs}$ as the linear system size, and $\text{nnz}.$ the number of non-zero values for the stiffness matrix. A few remarks: (i) $|\Theta|$ is smaller than $\text{nnz}.$ (by a factor of $10$ for $\lambda_{56}$), (ii) the metric error is higher for PINNs by more than a factor of $10$, and (iii) FEM keeps outperforming PINNs. Concerning execution time, LU and PINNs show more similar results, and GMRES is the best solution, as direct inversion does not scale well with dimension.

\begin{table}[!htb]
\renewcommand\arraystretch{1.3}
\begin{center}
\footnotesize
%\resizebox{8.75cm}{!} {
\begin{tabular}{
    |>{\centering\arraybackslash}m{2cm}% instead of "p" is "m"
    |>{\centering\arraybackslash}m{1cm}
    |>{\centering\arraybackslash}m{1cm}
    |>{\centering\arraybackslash}m{1.4cm}
    |>{\centering\arraybackslash}m{1.4cm}
    |>{\centering\arraybackslash}m{1.4cm}|
    }
    \hline    
\multicolumn{2}{|c|}{Setting} & $|\mT|$  & $|\Theta|$ &  Metric & $t_\text{exec}$ (s)\\ \hline
\multirow{2}{*}{PINN} & $\lambda_{92}$&  \multirow{2}{*}{8.0$\times$10$^3$} &  3.9$\times$10$^5$ &  3.7$\times$10$^{-1}$ & 2.4$\times$10$^3$ \\ \cline{2-2}\cline{5-6}
& $\lambda_{56}$ & & 1.7$\times$10$^5$ & 8.3$\times$10$^{-2}$& 9.3$\times$10$^2$ \\ \hline\hline
\multicolumn{2}{|c|}{Setting} & $N_\text{dofs}$ & nnz. & Metric & $t_\text{exec}$ (s)\\ \hline
\multirow{2}{*}{FEM} & LU  &   \multirow{2}{*}{8.0$\times$10$^3$} & \multirow{2}{*}{1.9$\times$10$^6$}  & 5.7$\times$10$^{-3}$ & 1.2 $\times$10$^2$\\  \cline{2-2}\cline{5-6}
& GMRES &&& 5.9$\times$10$^{-3}$  & 3.9 \\ \hline
 \end{tabular}%}
\end{center}\caption{Training v/s assembly+solver times : PINNs (left) v/s FEM (right).}
\label{tab:3DFEM} 
\end{table}

\section{Conclusion}\label{sec:Conclusion}
In this work, we  applied HPO via Gaussian processes-based Bayesian optimization to enhance the training of PINNs. We focused on forward problems for the Helmholtz operator, and carried out complete numerical experiments, with respect to performance, $(h,\kappa)$-analysis, and dimension. We compared the fitted PINNs with FEM. Numerical results: (i) confirm the performance and necessity of HPO, (ii) give a further insight on which could be good HPs for this problem, and (iii) pave the way toward running more efficient PINNs.

Further research include application to other operators, and enhancing the presented tuning procedure. Also, notice the active research area on PINNs theory. For example, error bounds \cite{Mishra2020EstimatesOT} would be valuable to obtain optimality conditions for $(h,\kappa)$ and to obtain bounded generalization and training errors. Both theoretical and empirical perspectives will allow to get closer to that objective. Also, it would be interesting to perform $\kappa$-analysis \cite{diwan2019can,spence2014}, and to introduce compression procedure for HPO.

\section*{Acknowledgement}
The authors would like to thank FES-UAI postdoc grant, ANID PIA/BASAL FB0002, and ANID/PIA/\newline ANILLOS ACT210096, for financially supporting this research.

\bibliography{references}


%\begin{thebibliography}{00}


%\bibitem{}

%\end{thebibliography}
\end{document}